\documentclass{elsart3-1}

\usepackage{amsmath,amssymb,graphicx,epsfig,color}
\usepackage[english,francais]{babel}
\usepackage{enumerate}
\newtheorem{theorem}{Theorem}%[section]

\newtheorem{e-proposition}[theorem]{Proposition}

\newtheorem{e-definition}[theorem]{Definition\rm}
\newtheorem{remark}{\it Remark\/}

\renewcommand{\theequation}{\arabic{equation}}
\def\rr{\mathbb{R}}

\def\l{\lambda}
\def\d{\partial}
\def\D{\Delta}
\def\q{\qquad}

\def\n{\nabla}
\def\og{\leavevmode\raise.3ex\hbox{$\scriptscriptstyle\langle\!\langle$~}}
\def\fg{\leavevmode\raise.3ex\hbox{~$\!\scriptscriptstyle\,\rangle\!\rangle$}}

\def\hoi{H_0^1(\Omega)}
\def\into{\int_{\Omega}}

\def\eps{\epsilon}

\begin{document}

\begin{frontmatter}

\selectlanguage{english}
\title{Hardy inequality and Pohozaev identity
for operators with boundary singularities:
 some applications}

\vspace{-2.6cm} \selectlanguage{francais}
\title{L'in\'{e}galit\'{e} de Hardy et l'identit\'{e} de Pohozaev pour op\'{e}rateurs avec
des singularit\'{e}s sur la fronti\`{e}re: quelques applications}

\selectlanguage{english}
\author[BCAM,UAM]{Cristian Cazacu},
\ead{cazacu@bcamath.org}

\address[BCAM]{BCAM - Basque Center for Applied Mathematics,
Bizkaia Technology Park, Building 500, 48160 Derio, Spain}
\address[UAM]{Departamento de Matem\'{a}ticas, Universidad Aut\'{o}noma,
 28049 Madrid, Spain}
\begin{abstract}
 We consider the Schr\"{o}dinger operator $A_\l:=-\D -\l/|x|^2$, $\l\in \rr$,
 when the singularity is located on the boundary
 of a smooth domain $\Omega\subset \rr^N$, $N\geq 1$

The aim  of this Note is two folded. Firstly, we justify the
extension of the classical Pohozaev identity for the Laplacian to
this case.  The problem we address is very much related to
Hardy-Poincar\'{e}
inequalities with boundary singularities. %As a consequence we show
%non-existence results for solutions of non-linear singular PDEs.
Secondly, the new Pohozaev identity allows to develop the multiplier
method for the wave and the Schr\"{o}dinger equations. In this way
we extend to the case of boundary singularities well known
observability and control properties for the classical wave and
Schr\"{o}dinger equations when the singularity is placed in the
interior of the domain (Vanconstenoble and Zuazua \cite{judith}).

 \vskip 0.5\baselineskip
\selectlanguage{francais} \noindent{\bf R\'esum\'e
 } \vskip 0.5\baselineskip \noindent

Nous allons consid\'{e}rer l'operateur de Schr\"{o}dinger $A_\l:=-\D
-\l/|x|^2$, $\l\in \rr$, quand l'origine est situ\'{e}e sur la
fronti\'{e}re d'un domaine born\'{e} et r\'{e}guli\`{e}re
$\Omega\subset \rr^N$, $N\geq 1$.

 Cette Note a deux objectifs. Premi\`{e}rement,
nous montrons la extension de l'identit\'{e} classique de Pohozaev
pour le Laplacien dans ce cas.
 Le
probl\`{e}me que nous abordons est tr\`{e}s li\'{e} aux
in\'{e}galit\'{e}s de Hardy-Poincar\'{e} avec des singularit\'{e}s
sur la fronti\'{e}re. %Comme cons\'{e}quence, nous d\'{e}montrons des
%r\'{e}sultats de non-existence pour les solutions des EDPs
%non-lin\'{e}aires et singuli\`{e}res.
En second lieu, la nouvelle identit\'{e} de Pohozaev permet
d\'{e}river le methode de multiplicateurs  pour les \'{e}quations
des ondes et de Schr\"{o}dinger. De cette fa\c{c}on, nous
\'{e}tendons au cas de la singularit\'{e}  frontali\`{e}re
propri\'{e}t\'{e}s d'observabilit\'{e} et contr\^{o}le  pour les
classiques \'{e}quations des ondes et de Schr\"{o}dinger bien connus
dans le cas d'une singularit\'{e} \`{a} l'interieur (Vancostenoble
et Zuazua \cite{judith})
\end{abstract}
\end{frontmatter}

% now the Version fran?aise abr?g?e, if it exists
\selectlanguage{francais}

\noindent {\bf Version fran\c{c}aise abr\'eg\'ee}\\

\begin{figure}[h!]
\begin{center}
 %$\boxed{\mathbf{\gamma_{\textrm{\bf loc}}=\gamma_{\textrm{\bf global}}=\gamma>0}}$
\setlength{\unitlength}{0.4cm}
\begin{picture}(16,7)
\linethickness{0.3mm}
  \put(2,3){\vector(1,0){11}}
  \put(7,1.3){\vector(0,1){7.2}}
 % \textcolor{blue}{\qbezier(3.2,6.8)(7.25,-1)(11,7)}
   \textcolor{red}{\qbezier(5,4)(7,1.9)(9,4.3)
   \qbezier(5,4)(2.9,6.5)(6.5,7.1)
   \qbezier(6.5,7.1)(11.9,8.6)(9,4.3)}
   \put(7,9.4){\makebox(0,0){\emph{\bf C1: The elliptic case }}}
   \put(7.4,4.5){\makebox(0,0){\emph{$\mathbf{\Omega}$}}}
   \put(11.8,2.5){\makebox(0,0){$\mathbf{x_N=0}$}}
   \put(5.2,8.2){\makebox(0,0){$\mathbf{x'=\bf{0}}$}}
  % \put(3.1,4.8){\makebox(0,0){$\mathbf{P_\gamma}$}}
  % \put(12.7,5.5){\makebox(0,0){$\mathbf{x_N=\gamma|x'|^2}$}}
    \put(6.3,2.2){\bf{0}}
%    \put(6,0.7){\bf{Figure 1}}
    \end{picture}
\begin{picture}(14,10)
\linethickness{0.3mm}
  \put(2,3){\vector(1,0){12}}
  \put(7.15,1){\vector(0,1){7.8}}
  %\textcolor[rgb]{0.25,0.00,0.25}{\qbezier(3,7)(7.3,-1.1)(11,7)}
  \textcolor{red}{
  \qbezier(5,4.1)(7.2,1.8)(9,4.3)
   \qbezier(5,4.1)(4,5)(3,1.5)
   \qbezier(3,1.5)(2.2,0)(4,5)
   \qbezier(11,1.5)(11.9,0)(10,5)
    \qbezier(4,5)(7,11)(10,5)
   \qbezier(9,4.3)(10,5.3)(11,1.5)}
  %\textcolor{blue}{ \qbezier(1.9,0.5)(7,5.55)(11.8,0.3)}
 %  \put(1.2,0.1){$\mathbf{P_{\gamma_{\textrm{\bf global}}}}$}
  % \put(9.7,0.0){$\mathbf{x_N=\gamma_{\textrm{\bf global}}|x'|^2}$}
  % \qbezier(5,4)(3,7)(6.5,5)
  % \qbezier(6.5,5)(11.5,9)(9,4.3)
   \put(7,9.5){\makebox(0,0){\emph{\bf C2: The locally elliptic case}}}
   \put(7.4,4.5){\makebox(0,0){\emph{$\mathbf{\Omega}$}}}
   \put(6.5,2.2){\bf{0}}
   \put(12.3,2.5){\makebox(0,0){$\mathbf{x_N=0}$}}
   \put(5.1,8.4){\makebox(0,0){$\mathbf{x'=\bf{0}}$}}
  % \put(1.6,6.8){\makebox(0,0){$\mathbf{P_{\gamma_{\textrm{\bf loc}}}}$}}
  % \put(13.3,6.5){\makebox(0,0){$\mathbf{x_N=\gamma_{\textrm{\bf loc}}|x'|^2}$}}
  % \put(6,0.3){\bf Figure 3}
  \end{picture}
\begin{picture}(15,10)
\linethickness{0.3mm}
  \put(2,3){\vector(1,0){11}}
  \put(7,1){\vector(0,1){7}}
%  \textcolor{blue}{\qbezier(2.6,1.2)(6.9,4.64)(11.5,1.4)}
   \textcolor{red}{\qbezier(4,2.4)(7,3.65)(10,2.3)
   \qbezier(4,2.4)(3.15,1.9)(5.1,5)
   \qbezier(10,2.3)(11,2)(9,5)
   \qbezier(5.1,5)(7.1,7.5)(9,5)}
  % \put(1.8,1.8){$\mathbf{P_\gamma}$}
   \put(6.5,2.1){\bf{0}}
  % \put(10.1,0.8){$\mathbf{x_N=\gamma|x'|^2}$}
   \put(7,9){\makebox(0,0){\emph{\bf C3: The hyperbolic case}}}
   \put(7.4,4.5){\makebox(0,0){\emph{$\mathbf{\Omega}$}}}
   \put(12.1,2.6){\makebox(0,0){$\mathbf{x_N=0}$}}
   \put(5.3,7.7){\makebox(0,0){$\mathbf{x'=\bf{0}}$}}
  \end{picture}
  \begin{picture}(14,10)
\linethickness{0.3mm}
  \put(2,3){\vector(1,0){12}}
  \put(7.15,1){\vector(0,1){7.8}}
  %\textcolor[rgb]{0.25,0.00,0.25}{\qbezier(3,7)(7.3,-1.1)(11,7)}
  \textcolor{red}{
 \qbezier(7,3)(11.4, 2.9)(10, 5)
   \qbezier(5,1.9)(6,3.1)(7,3)
   \qbezier(5,1.9)(2.2,0)(4,5)
   %\qbezier(11,1.5)(11.9,0)(10,5)
    \qbezier(4,5)(7,11)(10,5)
   %\qbezier(9,4.3)(10,5.3)(11,1.5)
   }
  %\textcolor{blue}{ \qbezier(1.9,0.5)(7,5.55)(11.8,0.3)}
 %  \put(1.2,0.1){$\mathbf{P_{\gamma_{\textrm{\bf global}}}}$}
  % \put(9.7,0.0){$\mathbf{x_N=\gamma_{\textrm{\bf global}}|x'|^2}$}
  % \qbezier(5,4)(3,7)(6.5,5)
  % \qbezier(6.5,5)(11.5,9)(9,4.3)
   \put(7,9.5){\makebox(0,0){\emph{\bf C4: $x\cdot \nu$ changes sign at 0}}}
   \put(7.4,4.5){\makebox(0,0){\emph{$\mathbf{\Omega}$}}}
   \put(6.6,2.2){\bf{0}}
   \put(12.3,2.5){\makebox(0,0){$\mathbf{x_N=0}$}}
   \put(5.1,8.4){\makebox(0,0){$\mathbf{x'=\bf{0}}$}}
  % \put(1.6,6.8){\makebox(0,0){$\mathbf{P_{\gamma_{\textrm{\bf loc}}}}$}}
  % \put(13.3,6.5){\makebox(0,0){$\mathbf{x_N=\gamma_{\textrm{\bf loc}}|x'|^2}$}}
  % \put(6,0.3){\bf Figure 3}
  \end{picture}\quad
\caption{\label{fig3} Geometry of $\Omega$ (in red), $\nu$ denotes the outward normal vector} %The case of local ellipticity
%$\gamma_{\textrm{loc}}>0$}
\end{center}
\end{figure}
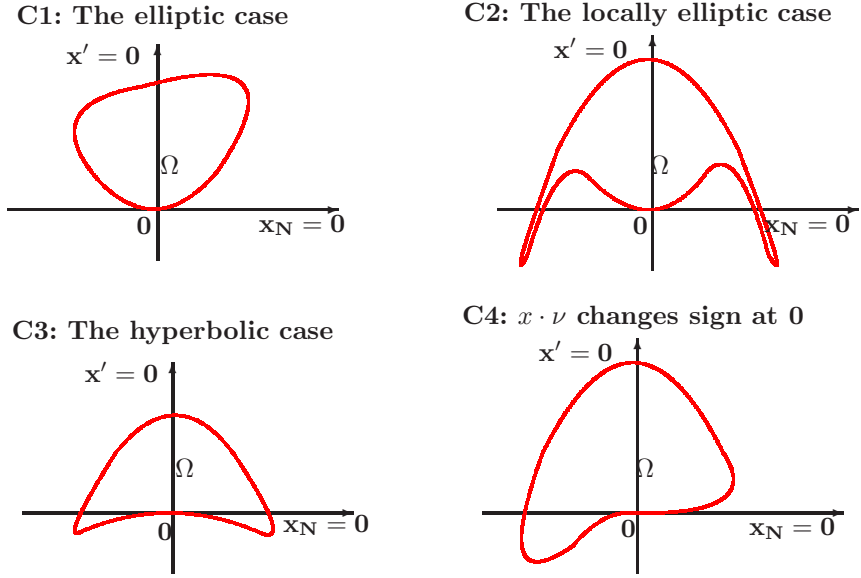

Dans cette Note nous nous int\'{e}ressons \`{a} l'op\'{e}rateur
$A_\l:=-\D-\l/|x|^2$, $\l \in \rr$, lorsque l'origine est situ\'{e}e
sur la fronti\`{e}re d'un domaine r\'{e}gulier $\Omega\subset
\rr^N$, $N\geq 1$. Il est connu que la valeur $\l(N):=N^2/4$, qui
est la meilleure constante dans les in\'{e}galit\'{e}s de Hardy
ci-apr\`{e}s, est critique lorsque l'on \'{e}tudie les
propri\'{e}t\'{e}s qualitatives de $A_\l$. Dans la premi\`{e}re
partie de cette Note nous montrons que pour tout $\l\leq \l(N)$,
l'identit\'{e} de Pohozaev (voir Th\'{e}or\`{e}me \ref{poho}) a lieu
dans le domaine de $A_\l$, d\'{e}fini par $D(A_\l):=\{ u\in H_\l\ |\
A_\l u \in L^2(\Omega)\}$. Nous allons d\'{e}finir plus tard
l'espace $H_\l$ et quelques unes de ses propri\'{e}t\'{e}s.
Formellement, le Th\'{e}or\`{e}me \ref{poho} peut \^{e}tre obtenu
par int\'{e}gration directe. Cependant, la singularit\'{e} $x=0$
engendre une perte de r\'{e}gularit\'{e} de l'op\'{e}rateur $(A_\l,
D(A_\l))$ et les int\'{e}grations par parties ne sont plus
justifi\'{e}es rigoureusement. De plus, la r\'{e}gularit\'{e} $L^2$
de la d\'{e}riv\'{e}e normale n'a plus lieu car les estimations
elliptiques standards ne s'appliquent plus puisque la
singularit\'{e} est localis\'{e}e sur le bord. N\'{e}anmoins, la
trace d'un \'{e}l\'{e}ment de $(A_\l, D(A_\l))$ existe dans un
espace $L^2$ \`{a} poids, dont le poids est g\'{e}n\'{e}r\'{e} \`{a}
l'origine, comme il est montr\'{e} dans le Th\'{e}or\`{e}me
\ref{trace}. Dans la deuxi\`{e}me partie de cette Note, nous
montrons plusieurs applications des Th\'{e}or\`{e}mes \ref{trace},
\ref{poho}. D'abord des solutions non-triviales d'une EDP
singuli\`{e}re sont trait\'{e}es dans le Th\'{e}or\`{e}me
\ref{teo4}. Ensuite, nous d\'{e}rivons des techniques de
multiplicateurs afin de prouver "la r\'{e}gularit\'{e} cach\'{e}e"
de la d\'{e}riv\'{e}e normale dans le cas de l'\'{e}quation des
ondes et de l'\'{e}quation de Schr\"{o}dinger, correspondant
\`{a}$A_\l$ (voir le Th\'{e}or\`{e}me \ref{Wmult}). En particulier,
nous r\'{e}pondons \`{a} la question concernant la
controlabilit\'{e} des syst\`{e}mes conservatifs. Le r\'{e}sultat
principal est donn\'{e} par le Th\'{e}or\`{e}me \ref{t1} et est
d\^{u} \`{a} l'identit\'{e} des multiplicateurs \eqref{multipliers},
en combinaison avec une in\'{e}galit\'{e} forte de Hardy,
formul\'{e}e dans le Th\'{e}or\`{e}me \ref{tu8}. Pour plus de
claret\'{e} dand la pr\'{e}sentation, nous allons traiter notamment
le cas C1 de la figure \ref{fig3}. Cependant, les m\^{e}mes
r\'{e}sultats peuvent \^{e}tre \'{e}tendus aux cas C2, C3, C4 dans
un cadre fonctionnel plus faible, d\^{u} \`{a} l'in\'{e}galit\'{e}
plus faible de Hardy \eqref{oeq44}.

\smallskip
 \selectlanguage{english}
\noindent {\bf 1.  Introduction}\\
 Let us consider
$\Omega$ to be a smooth subset of $\rr^N$, $N\geq 1$, with the
origin $x=0$ placed on its boundary $\Gamma$. Without losing the
generality we distinguish  the following geometrical configurations
for $\Omega$ as in Figure \ref{fig3}:  \noindent {\bf C1 -} $\Omega$
is a subset of $\rr_{+}^{N}:=\{x=(x_1,\ldots, x_N)\in \rr^N \ |\
x_N>0 \}$ (Fig. \ref{fig3}, top left). {\bf C2 -} Close to $x=0$,
the points $x\in \Gamma$ satisfy $x\cdot \nu \geq 0$. Nevertheless,
$\Omega$ crosses the hyperplane
 $x_N=0$ far from origin (Fig. \ref{fig3}, top, right).
 {\bf C3 -} Close to $x=0$, the points
  $x\in \Gamma$ verify $x\cdot \nu\leq 0$ (Fig. \ref{fig3}, bottom, left). {\bf
C4 -}  For $x\in \Gamma$ the sign of $x\cdot \nu$ changes at origin
(Fig. \ref{fig3}, bottom, right).

The following Hardy inequalities are well-known: if $\Omega$
verifies the case C1 in Fig. \ref{fig3}, then (e.g.
\cite{cristiCRAS}) for any $u\in C_{0}^{\infty}(\Omega)$ it holds
that
\begin{align}\label{oeq3}
\into |\nabla u|^2dx\geq \frac{N^2}{4}\into
\frac{u^2}{|x|^2}dx+\frac{1}{4}\into
\frac{u^2}{|x|^2\log^2(R_{\Omega}/|x|)}dx.
\end{align}
 where $R_\Omega=\sup_{x\in \overline{\Omega}}|x|$. If $\Omega$ satisfies the cases  C2, C3, C4 as in Fig.
\ref{fig3} then (e.g. \cite{Fall1}, \cite{musina}) there exist
$C_2=C_2(\Omega)\in \rr$ and $C_3=C_3(\Omega, N)>0$ such that  any
$u\in C_{0}^{\infty}(\Omega)$ satisfies
\begin{align}\label{oeq44}
C_2\into u^2dx +\into |\nabla u|^2dx\geq \frac{N^2}{4}\into
\frac{u^2}{|x|^2}dx +C_3\into
\frac{u^2}{|x|^2\log^2(R_\Omega/|x|)}dx.
\end{align}
In both situations above, the constant $\l(N)=N^2/4$ is optimal.
 In the sequel we consider that $\Omega$ verifies the case C1 in
Fig. \ref{fig3} and let be the operator $(A_\l, D(A_\l))$, $\l\leq
\l(N)$, acting on $\Omega$.

Firstly, this Note is aimed to justify  the  Pohozaev identity
\eqref{pohoidentity} in the functional setting $(A_\l, D(A_\l))$ in
which  Theorem \ref{trace} plays a crucial role. Pohozaev-type
identities (e.g. pp. 515, \cite{evans}) have been widely used to
show non-existence results  to nonlinear elliptic equations. In
particular, we point out that this issue has been also studied for
nonlinear equations with singular potentials (see e.g. \cite{Peral},
\cite{gosu}). In those cases, due to the regularizing effect of the
nonlinearity, the solutions become regular enough to obtain the
corresponding Pohozaev identity by direct computations. This is not
precisely our case. To the best of our knowledge, the regularity of
the operator $(A_\l, D(A_\l))$ in terms of the Sobolev spaces is yet
unknown and this affects the direct justification of
\eqref{pohoidentity}. \
To do this we proceed by approximations arguments near the
singularity. We analyze two situations when discussing the values of
$\l$: the subcritical case $\l\leq \l(N)$ respectively the critical
case $\l=\l(N)$. However,
 the main novelty appears for the critical value $\l=\l(N)$
 in which case $H_\l$ is strictly larger than $\hoi$. This case
 requires e better understanding of the norm  of $H_{\l(N)}$ as
 discussed in Section 2.
% Then the identity \eqref{pohoidentity}   allows to
% derive the multiplier method which is one of the main tools to be used
% in applications to control.

 The controllability properties and stabilization of the heat and
wave equation corresponding to $A_\l$ have been analyzed in
\cite{heatjudith}, \cite{sylvain}, \cite{judith}
 in the case of interior singularity. Roughly speaking, they showed
 that the parameter $\lambda_\star=(N-2)^2/4$, which is the optimal constant in the Hardy
inequality with interior singularity, is critical when asking the
well-posedness and control properties of such systems. In the second
part of this  Note we address the question of controllability for
the wave
and Schrodinger equations corresponding to $A_\l$, in the case of boundary singularity. %This
 Our main result asserts that  we can increase
the range of values $\lambda$ (from $\l_\star$ to $\l(N)$) for which
the exact controllability holds.

\smallskip
\noindent {\bf 2. The space $H_\l$,  proper norm and main elliptic
results.}
% Assume $\Omega\subset \rr^N$, $N\geq 1$ satisfies
%  the case C1 in Fig. \ref{fig3} and let $\l \leq \l(N)$.
Following the idea in  \cite{vazzua},  thanks to inequality
\eqref{oeq3} we consider the Hardy functional
%\begin{equation}\label{normB}
$B_{\l}[u]=\into \big[|\n u|^2 -\l u^2/|x|^2\big]dx$,
%\end{equation}
which is positive and finite for all $u\in C_{0}^{\infty}(\Omega)$.
We define the Hilbert space $H_\l$ to be  the completion of
$C_0^\infty(\Omega)$ functions in the norm induced by $B_\l[u]$.
 If
$\lambda<\lambda(N)$, it holds that
 $\hoi=H_\l$ due to Hardy inequality %according to the estimates
%$ \big(1-\lambda/\l(N)\big) ||u||_{\hoi}^2 \leq B_\l[u]\leq
%||u||_{\hoi}^{2}$,
which ensures the equivalence of the norms. Similar to the case of
interior singularity emphasized in \cite{hiddenenergy}, an
interesting phenomena appears  in the critical space $H_{\l(N)}$.
%The functional $B_{\l(N)}[u]$ makes sense as an improper integral
%for more singular distributions.
 Assume $\Omega=\{x\in \rr^N |\ |x|\leq 1, \ x_N>0\}$ and let $e_1:=x_N|x|^{-N/2}J(z_{0,1}|x|)$ where
$z_{0,1}$ is the first positive zero of the Bessel function $J_0$.
 Then $\textrm{there exists }\lim_{\eps\rightarrow 0}\int_{x\in \Omega, |x|\geq
\eps}\big[ |\n e_1|^2-\l(N)e_1^2|x|^2\big]dx< \infty$ although
$e_1\not \in \hoi$. Surprisingly, the meaning of $B_{\l(N)}[\cdot]$
is not well-defined in the sense of principle value. Indeed, if it
were one can check that $B_{\l(N)}[e_1-\phi]\geq C_0>0$, for all
$\phi \in C_{0}^{\infty}(\Omega)$ for some universal constant
$C_0>0$.
This is in contradiction with the definition of $H_{\l(N)}$ ! %Next
The remedy for this is to consider the functional
\begin{equation}\label{norm2SB}
B_{\l, 1}[u]=\int_\Omega \Big|\n u+ \frac{N}{2}
\frac{x}{|x|^2}u-\frac{e_N}{x_N}u\Big|^2dx+ (\l(N)-\l)\into
\frac{u^2}{|x|^2}dx,
\end{equation}
where $e_N$ denotes the $N-$th canonical vector of $\rr^N$. We
observe that $B_\l[u]=B_{\l, 1}[u], \textrm{ for all } u\in
C_{0}^{\infty}(\Omega).$
  Therefore, $B_{\l, 1}[u]$ induces a new norm in the space $H_\l$
   which is well
understood in the sense of principal value.
 In the sequel, we denote
by
$||\cdot||_{H_\l}$ the norm induced by $B_{\l, 1}$.\\
\noindent {\bf Notations:} % The domain of $A_\l$ is $D(A_\l):=\{ u\in H_\l \
%|\ A_\l u \in L^2(\Omega)\}.$
For any $\eps>0$, $\theta_\eps$ is a smooth cut-off function which
satisfies $\theta_\eps=0$ for $|x|\leq \eps$ respectively
$\theta_\eps=1$ for $|x|\geq 2\eps$. Besides, $\vec{q}\in
(C^2(\overline{\Omega}))^N$ denotes a vector field such that
$\vec{q}=\nu$ on $\Gamma$.  $H_{\l}^{'}$ denotes the dual space of
$H_\l$.  Next we state our main elliptic results.
\smallskip
\begin{theorem}[trace regularity]\label{trace}
Assume $\Omega\subset\rr^N$, $N\geq 1$, verifies the case C1. Let us
consider $\l \leq \l(N)$ and $u\in (A_\l, D(A_\l))$. Then
%\begin{equation}\label{tracereg}
$\big(\d u/\d \nu\big)|x|\in L^2(\Gamma)$
%\end{equation}
and there exists  $C=C(\Omega)>0$ such that
\begin{equation}\label{ineqtrace}
\int_{\Gamma}\Big(\frac{\d u }{\d \nu}\Big)^2|x|^2 d\sigma \leq
C(||u||_{H_\l}^2+||A_\l u||_{L^2(\Omega)}^2).
\end{equation}
\end{theorem}
\textsl{ Sketch of the proof.} In order to avoid the singularity, we
multiply  $A_\l u$ by  $|x|^2(\vec{q}\cdot \n u) \theta_\eps$
%  (such an
%election of $\vec{q}$ can be always done in smooth domains, see
%\cite{MR953547}, Lemma 3.1, page 29).
and we integrate by parts. Then we obtain an identity which,
combined with Cauchy-Schwartz inequality, allows to get uniform
upper bounds for the boundary term. Then, by Fatou Lemma we can pass
to the limit as $\eps$ tends to zero to end up the proof. $\square$
\smallskip
\begin{theorem}[Pohozaev identity]\label{poho}
Assume $\Omega\subset \rr^N$, $N\geq 1$, verifies the case C1 and
let $\l\leq \l(N)$. Then for all $u\in (A_\l, D(A_\l))$ it holds
that
\begin{align}\label{pohoidentity}
\frac{1}{2}\int_{\Gamma} (x\cdot \nu)\Big(\frac{\d u}{\d \nu}\Big)^2
d\sigma&= -\into A_\l u(x\cdot \n u)dx- \frac{N-2}{2}||u||_{H_\l}^2.
\end{align}
\end{theorem}
\textsl{Sketch of the proof.} Note firstly that all terms in
 \eqref{pohoidentity} are finite. Indeed, thanks to the fact that $x\cdot
 \nu=O(|x|^2)$ and Theorem \ref{trace} we obtain the integrability of
 the boundary term. Moreover, $x\cdot \n u\in L^2(\Omega)$ as shown in
 Theorem \ref{tu8}.
  We proceed by multiplier technique. If $\l<
\l(N)$ the multiplier $(x\cdot \n u)\theta_\eps$ is used to obtain
\eqref{pohoidentity} in the limit process. If $\l=\l(N)$  we apply
the subcritical result for $(\l(N)-\delta)$ and we argue we can pass
to the limit as $\delta$ tends to zero. Indeed, following
approximation lemma which is proved in generality in \cite{acz}, for
a given $u\in D(A_\l)$, the solution $u_\delta$ of
$A_{\l(N)-\delta}u_\delta=A_{\l(N)}u$, converges strongly to $u$ in
$H_{\l(N)}$. Secondly, by comparison arguments we show the
asymptotic behavior of $u_\delta$, $\n u_\delta$ when $A_\l u$ is
smooth. This is done by constructing super solutions and using
rescalling arguments in \cite{marcus2} to get rid of the
singularity. A density argument together with Theorem \ref{trace}
conclude the proof. $\square$
\smallskip
\begin{theorem}\label{teo4}
 Assume $\l\leq  \l(N)$ and %let us consider
$\Omega\subset \rr^N$, $N\geq 3$, %be a star-shaped domain (i.e. $x\cdot \nu\geq
%0$, for all $x\in \Gamma$)
 satisfying the case C1.  Let us consider the problem
\begin{equation}\label{nonex}  -\D u-\frac{\lambda}{|x|^2}u=|u|^{\alpha-1}u,  \quad  x\in \Omega; \q
  u=0, \quad  x\in \Gamma,
  \end{equation}
1). If $1<\alpha< (N+2)/(N-2)$,  problem \eqref{nonex} has non
trivial solutions in $H_\l$. Moreover, if $1< \alpha< N/(N-2)$
it has non trivial solutions in $(A_\l, D(A_\l))$. \\
2). Assume $\Omega$ is a star-shaped domain (i.e. $x\cdot \nu\geq
0$, for all $x\in \Gamma$). If $\alpha \geq (N+2)/(N-2)$,  problem
\eqref{nonex} does not have non trivial solutions in $(A_{\l},
D(A_{\l}))$.
\end{theorem}
\textsl{Sketch of the proof.} The existence of non-trivial solutions
reduces to finding minimizers for the problem
$I=\inf_{||u||_{L^{\alpha+1}(\Omega)}=1}||u||_{H_\l}^{2}$ and is due
to the compact embedding $H_\l\hookrightarrow L^{\alpha+1}(\Omega)$.
The non-existence result is proved by combining Theorem \ref{poho}
and unique continuation results as in \cite{kenig}. $\square$.

\smallskip
 \noindent {\bf 3.  Applications to controllability}\label{5sec}

%
%In the sequel we assume $\Omega\subset \rr^N$, $N\geq 1$, satisfies
%case C1 in Fig. \ref{fig3}. Moreover,
Let us  set
%\begin{equation}\label{oeq2}
$\Gamma_0:=\{x\in \Gamma\ |\ x \cdot \nu \geq 0\}.$
%\end{equation}
 a non-empty part of the boundary $\Gamma$. Next we consider the Wave-like process
\begin{equation}\label{eq124}\left\{\begin{array}{ll}
  u_{tt}-\Delta u-\lambda \frac{u}{|x|^2}=0, & (t,x)\in Q_T, \\
  u(t,x)=h(t,x)\chi_{\Gamma_0}, & (t,x)\in \Sigma_T, \\
  (u(0,x), u_t(0,x))=(u_0(x), u_1(x)), &  x\in \Omega, \\
\end{array}\right.
\end{equation}
where  $Q_T=(0, T)\times \Omega$, $\Sigma_T=(0,T)\times \Gamma$ and
$\chi_{\Gamma_0}$ denotes the characteristic function of $\Gamma_0$.

 The solution of (\ref{eq124}) is defined in week sense by the transposition
 method (J.L. Lions \cite{lions1}).
 In this section we address the question of  exact controllability of system (\ref{eq124}),
 i.e. whether for any initial data
$(u_0,u_1)\in L^2(\Omega)\times H_{\l}^{'}$  and any target
$(\overline{u_0},\overline{u_1})\in L^2(\Omega)\times H_{\l}^{'}$ ,
there exists a finite time $T>0$ and a control $h\in L^2((0,T)\times
\Gamma_0)$ such that the solution of (\ref{eq124}) satisfies
$(u_t(T,x),u(T,x))=(\overline{u_1}(x),\overline{u_0}(x))$ for all
$x\in \Omega$.  In view of the time-reversibility of the equation it
is enough to consider the null-controllability problem, i.e. the
case where the target $(\overline{u_0},\overline{u_1})=(0,0)$.
 By now classical \textit{Hilbert Uniqueness Method} (HUM) (see J. L. Lions \cite{lions1}) the null-controllability of system (\ref{eq124}) is
characterized through the adjoint system
\begin{equation}\label{eq2}\left\{\begin{array}{ll}
  v_{tt}-\Delta v-\lambda \frac{v}{|x|^2}=0, & (t,x)\in Q_T, \\
  v(t,x)=0, & (t,x)\in \Sigma_T, \\
  (v(0,x), v_t(0,x))=(v_0(x), v_1(x)), &  x\in \Omega, \\
\end{array}\right.
\end{equation}
The operator $(\mathcal{A_\l}, D(\mathcal{A_\l}))$ defined by
$\mathcal{A_\l}(v,w) = (w, \D v + \l |x|^2v)$ for all $(v,w) \in
D(\mathcal{A_\l}) = D(A_\l) \times H_\l$ generates the wave
semigroup in $H_\l\times L^2(\Omega)$. In view of that,  the adjoint
system is well-posed.

%\subsubsection{Hidden regularity for the normal derivative. From Pohozaev-type identity to Multipliers identity.}
 % Proof of main result (Theorem \ref{t1})}
In the sequel, we justify some ``hidden regularity" effect for the
system \eqref{eq2} which may not be directly deduce from the
semigroup regularity but from the equation itself.
\smallskip
\begin{theorem}[\bf Hidden regularity]\label{Wmult}
Assume $\l \leq \l(N)$ and $v$ is the solution of \eqref{eq2}
corresponding to the initial data $(v_0, v_1)\in H_\l\times
L^2(\Omega)$. Then $v$ satisfies
\begin{equation}\label{ineqtrace1}
\int_{0}^{T}\int_{\Gamma}\Big(\frac{\d v }{\d \nu}\Big)^2|x|^2
d\sigma dt \leq C(||v_0||_{H_\l}^2+||v_1||_{L^2(\Omega)}^2).
\end{equation}
for some universal constant $C>0$. Moreover, $v$ verifies the
identity
\begin{equation}\label{multipliers}
\frac{1}{2}\int_{0}^{T}\int_{\Gamma} (x\cdot \nu)\Big(\frac{\d v
}{\d \nu}\Big)^2d\sigma dt =\frac{T}{2}
(||v_0||_{H_\l}^2+||v_1||_{L^2(\Omega)}^2) +\int_{\Omega} v_t
\big(x\cdot \n v+\frac{N-1}{2}v\big)\Big|_{0}^{T}dx.
\end{equation}
\end{theorem}
\textsl{ Sketch of the proof.} By density, tt suffices to  prove
Theorem \ref{Wmult}, for initial data $(v_0, v_1)$ in
$D(\mathcal{A_\l})=D(A_\l)\times H_\l$. %Then, by density, we extend
%the result for initial data $(v_0, v_1)\in H_\l \times L^2(\Omega)$.
For the proof of \eqref{ineqtrace1} we multiply $A_\l v$ by
$|x|^2(\vec{q}\cdot \n v) \theta_\eps$ and integrate. The
integration in time and the conservation of energy allow to obtain
uniform bounds for the boundary term  in the energy space. Then by
Fatou Lemma we pass to the limit as $\eps$ tends to zero and the
proof finishes.

For the proof of \eqref{multipliers} we proceed straightforward from
Theorem \ref{poho}. Indeed, %according to semigroup theory of
%evolution equations we have $v \in C([0,T]; D_\lambda)\cap
%C^1([0,T]; H_\l)\cap C^2([0,T]; L^2(\Omega)).$
  %Next,
   for a fixed time $t\in [0, T]$ we apply
    Theorem \ref{pohoidentity} for $A_\l v=-v_{tt}$. Then we
   integrate in time, and due to the equipartition of energy
     we can finish the proof. $\square$

Due to Theorem \ref{Wmult} the operator $(v_0, v_1)\mapsto
\big(\int_{0}^{T}\int_{\Gamma_0}(x\cdot \nu)(\d v/ \d \nu)^2d\sigma
dt\big)^{1/2}$ %is a norm in $D(A_\l)\times H_\l$ but can be extended
 is a linear continuous map in $H_\l\times L^2(\Omega)$. Let
$\mathcal{H}$ be the completion of this norm in $H_\l \times
L^2(\Omega)$.  We consider the functional $J: \mathcal{H}\rightarrow
\rr$ defined by
\begin{equation}\label{functional}
J(v_0, v_1)(v):=\frac{1}{2}\int_{0}^{T}\int_{\Gamma_0}(x\cdot
\nu)\Big(\frac{\partial v} {\partial \nu}\Big)^2d\sigma dt
-<u_1,v_0>_{H_{\lambda}^{'},H_\lambda}+(u_0,v_1)_{L^2(\Omega),L^2(\Omega)},
\end{equation}
where $v$ is the solution of \eqref{eq2} corresponding to initial
data $(v_0, v_1)$. Of course,
$<\cdot,\cdot>_{H_{\lambda}^{'},H_{\lambda}}$ denotes the duality
product.  The control $h\in L^2((0, T)\times \Gamma_0)$ for
\eqref{eq124}  could be chosen as $h=(x\cdot \nu) v_{\min}$ where
$v_{\min}$ minimizes the functional $J$ on $\mathcal{H}$ among the
solutions $v$ of (\ref{eq2}) corresponding to the initial data
$(u_0, u_1)\in H_{\lambda}^{'}\times L^2(\Omega)$ The existence of a
minimizer of $J$ %(and implicit the existence of a control $h$)
 is
assured by the coercivity of $J$ which is equivalent to
 the \textsl{Observability Inequality} for the adjoint
system  stated as follows.
\smallskip
\begin{theorem}[{\bf Observability}]\label{t1}
For all $\lambda\leq \lambda(N)$, there exists a positive constant
$D_1=D_1( \Omega, \l , T)>0$ such that for all $T\geq 2R_\Omega$,
and any initial data $(v_0, v_1)\in H_\l \times L^2(\Omega)$ the
solution of (\ref{eq2}) verifies
\begin{equation}\label{ObsIneq1}
%\int_{\Omega}\Big[ v_t^2(0,x)+|\nabla v(0,
%x)|^2-\l\frac{v^2(0,x)}{|x|^2}\Big]dx
||v_0||_{H_\l}^{2}+||v_1||_{L^2(\Omega)}^2\leq
D_1\int_{0}^{T}\int_{\Gamma_0}(x\cdot \nu)\Big(\frac{\d v}{\d
\nu}\Big)^2d\sigma dt.
\end{equation}
%For $\lambda=\lambda(N)$, there exists a positive constant $D_2(T,
%\Omega)$ such that for all $T\geq 2R_\Omega$, the solution of
%(\ref{eq2}) verifies
%\begin{equation}\label{ObsIneq2}
%%\int_{\Omega}\Big[ v_t^2(0,x)+|\nabla v(0,
%%x)|^2-\l\frac{v^2(0,x)}{|x|^2}\Big]dx
%E_{v}^{\l(N), 1}(0)\leq
%D_2(T,\Omega)\int_{0}^{T}\int_{\Gamma_0}(x\cdot \nu)\Big(\frac{\d
%v}{\d \nu}\Big)^2d\sigma dt.
%\end{equation}
\end{theorem}
\textsl{Sketch of the proof.}
 The proof of Theorem \ref{t1} relies
mainly on  Theorem \ref{Wmult},  combining compactness uniqueness
argument (cf. \cite{mach}) and the sharp Hardy inequality stated in
Theorem \ref{tu8}.
\smallskip
\begin{theorem}\label{tu8}
Assume $\Omega$  satisfies one of the cases C1-C4. Then, there
exists a constant $C=C(\Omega)\in \rr$ such that
\begin{equation}\label{equu41}
\int_{\Omega}|x|^2|\nabla v|^2dx\leq
R_{\Omega}^{2}\Big[\int_{\Omega}|\nabla
v|^2dx-\frac{N^2}{4}\int_{\Omega}\frac{v^2}{|x|^2}dx\Big]+C\int_{\Omega}v^2dx\quad
\forall v\in C_{0}^{\infty}(\Omega).
\end{equation}
\end{theorem}
\begin{remark}
The proof of Theorem \ref{tu8} is quite technical and we omit it
here. The constant $R_\Omega^2$ which appears in inequality
\eqref{equu41}, helps to obtain the control time $T>T_0=2 R_\Omega$,
which is sharp from the Geometric Control Condition considerations,
see \cite{BardosLeRa}.
\end{remark}
 The results above
guarantee the exact boundary controllability of (\ref{eq124}). More
precisely, we obtain
\smallskip
\begin{theorem}[{\bf Controllability}]\label{ht1}
Assume that $\Omega$ satisfies C1 and $\lambda\leq \lambda(N)$. For
any time $T>2R_{\Omega}$, $(u_0,u_1)\in L^2(\Omega)\times
H_{\lambda}^{'}$ and $(\overline{u_0}, \overline{u_1}) \in
L^2(\Omega)\times H_{\lambda}^{'}$ there exists $h\in
L^2((0,T)\times\Gamma_0)$ such that the solution of (\ref{eq124})
satisfies
%\begin{equation*}
$(u_t(T,x),u(T,x))=(\overline{u_1}(x),\overline{u_0}(x))$  for all
$x\in \Omega$.
%\end{equation*}
\end{theorem}
\smallskip
{\bf  Schr\"{o}dinger equation.}  In the above geometrical settings,
we consider the Schr\"{o}dinger equation
\begin{equation}\label{Seq124Sb}\left\{\begin{array}{ll}
  iu_{t}-\Delta u-\lambda \frac{u}{|x|^2}=0, & (t,x)\in Q_T, \\
  u(t,x)=h(t,x)\chi_{\Gamma_0}, & (t,x)\in \Sigma_T, \\
  u(0,x)=u_0(x), &  x\in \Omega, \\
\end{array}\right.
\end{equation}
%We define the Hilbert spaces $L^2(\Omega; \mathbb{C})$ and
%$H_0^1(\Omega; \mathbb{C})$ endowed with the inner products $<u,
%v>_{L^2(\Omega; \mathbb{C})}:=\textrm{Re}\into
%u(x)\overline{v(x)}dx$, $<u, v>_{H_0^1(\Omega;
%\mathbb{C})}:=\textrm{Re}\into \n u(x)\cdot \n \overline{v(x)} dx,$.
 For all $\l \leq
\l(N)$,  we  define the Hilbert space $H_\l (\Omega; \mathbb{C})$ as
the completion of $H_0^1 (\Omega; \mathbb{C})$ with respect to the
norm induced by the inner product $<u, v>_{H_\l(\Omega;
\mathbb{C})}:=\textrm{Re}\into \big(\n u(x)\cdot \n
\overline{v(x)}-\l u(x)\overline{v(x)}/|x|^2\big) dx$.  Then
\smallskip
\begin{theorem}[Controllability]
For any $\lambda\leq \l(N)$,  $u_0\in H_\l^{'}$,  $\overline{u_0}
\in H_{\lambda}^{'}$ and any time $T>0$ there exists  $h\in
L^2((0,T)\times\Gamma_0)$ such that the solution of \eqref{Seq124Sb}
 satisfies $u(T,x))= \overline{u_0}(x)$ for all $x\in \Omega$.
%\begin{remark}
\end{theorem}
  This result holds true due to the result valid for the wave equation.
Indeed, the general theory presented in an abstract form in
\cite{tuch},  assure the observability of systems like $\dot
z=iA_0z$ using results available for systems of the form
$\ddot{z}=-A_0 z$.
%\end{remark}

\smallskip\noindent\textbf{Acknowledgements.}
The author thanks Enrique Zuazua and Adi Adimurthi for fruitful
discussions and suggestions.

 Partially supported by the Grant MTM2008-03541 of the MICINN,
Spain, project
  PI2010-04 of the Basque Government, the ERC Advanced Grant FP7-246775
  NUMERIWAVES,
and a doctoral fellowship from UAM (Universidad Aut\'{o}noma de
Madrid).


\begin{thebibliography}{99}
%\bibitem{adimurthi1} N. C. Adimurthi and M. Ramaswamy, An improved Hardy-Sobolev
%inequality and its application, Proc. Amer. Math. Soc. 130 (2002),
%no. 2, 489-505 (electronic).
\bibitem{acz} A. Adimurthi, C. Cazacu and E. Zuazua, Best constants and Pohozaev
identity for Hardy-Sobolev type operators, in preparation.
\bibitem{BardosLeRa} C. Bardos,
 G. Lebeau and J. Rauch, Control
and stabilization for hyperbolic equations,
 Mathematical and numerical aspects of wave propagation
              phenomena ({S}trasbourg, 1991), 252-266,
 SIAM, Philadelphia, PA, 1991.
%\bibitem{MR1605678} H. Brezis and J. L. V\'{a}zquez, Blow-up solutions of some nonlinear
%elliptic problems, Rev. Mat. Univ. Complut. Madrid 10 (1997), no. 2,
%443-469.
\bibitem{marcus2} H. Brezis, M. Marcus and I. Shafrir, Extremal functions for {H}ardy's
inequality with weight, J. Funct. Anal., Journal of Functional
Analysis, 171, 2000,  1, 177-191.
\bibitem{cristiCRAS} C. Cazacu, On {H}ardy inequalities with
     singularities on the boundary, C. R. Acad. Sci. Paris, Ser. I, 349, 2011,
     273-277,
%\bibitem{simon} H. L. Cycon, R. G. Froese, W. Kirsch, and B. Simon, Schr\"{o}dinger
%operators with application to quantum mechanics and global geometry,
%Berlin (1987).
\bibitem{kenig}    D. Jerison and C. E. Kenig,
   Unique continuation and absence of positive eigenvalues for
              {S}chr\"odinger operators, With an appendix by E. M. Stein,
              Ann. of Math. (2), 121, 1985, 3, 463-494,
\bibitem{Peral} J. Davila and I. Peral,  Nonlinear elliptic problems with a singular weight
on the boundary, Calc. Var. doi 10.1007/s00526-010-0376-5.
\bibitem{sylvain} S. Ervedoza, Control and stabilization properties for a singular heat
              equation with an inverse-square potential,
              Comm. Partial Differential Equations, 33, 2008, 10-12,
1996--2019,
\bibitem{evans} L. C. Evans, Partial differential equations, Graduate
Studies in Mathematics, 19, American Mathematical Society,
Providence, RI, 2010, xxii+749.
\bibitem{Fall1} M. M. Fall, On the Hardy Poincar\'e inequality with boundary
 singularities. Commun. Contemp. Math, to appear.
\bibitem{musina} M. M. Fall and R. Musina,  Hardy-Poincar\'e inequality with
boundary singularities.  Proc. Roy. Soc. Edinburgh, to appear.
\bibitem{gosu} N. Ghoussoub, and C. Yuan, Multiple solutions for quasi-linear {PDE}s involving the
              critical {S}obolev and {H}ardy exponents, Trans. Amer. Math. Soc.,
    352, 2000, 12, 5703-5743.
\bibitem{hardy-polya} G. H. Hardy, J. E. Littlewood, and G. P\'{o}lya, Inequalities,
Cambridge Mathematical Library, Cambridge University Press,
Cambridge, 1988, Reprint of the 1952 edition.
%\bibitem{mouhamed} M. F. Mouhamed and R. Musina, Sharp nonexistence
%results for a linear elliptic inequality involving Hardy and Leray
%potentials, http://arxiv.org/abs/1006.5603.
%\bibitem{pot1} J.M. L\'{e}vy-Leblond, Electron capture by polar molecules,
%Phys. Rev., 153, (1) (1967) ,1-4.
\bibitem{lions1} J.-L. Lions, Contr\^olabilit\'e exacte,
 perturbations et stabilisation de
              syst\`emes distribu\'es. {T}ome 1,
   Recherches en Math\'ematiques Appliqu\'ees [Research in
              Applied Mathematics], 8, Contr{\^o}labilit{\'e} exacte.
               [Exact controllability], With appendices
               by E. Zuazua, C. Bardos, G. Lebeau and J.
              Rauch, Masson, Paris, 1988,
\bibitem{mach} E. Machtyngier, Exact controllability for
the {S}chr\"odinger equation, SIAM J. Control Optim., 32,
      1994, 1, 24--34,
%\bibitem{terti} A. Tertikas, S.  Filippas and J. Tidblom, On the structure of
%Hardy-Sobolev-Maz'ya inequalities, J. Eur. Math. Soc. (JEMS) 11
%(2009), no. 6, 1165-1185.
\bibitem{tuch} M. Tucsnak  and G. Weiss,
Observation and control for operator semigroups, Birkh\"auser
Advanced Texts: Basler Lehrb\"ucher.
              [Birkh\"auser Advanced Texts: Basel Textbooks],
              Birkh\"auser Verlag, Basel, 2009, xii+483.
\bibitem{heatjudith}, J. Vancostenoble, J. and E. Zuazua, Null controllability for the heat equation with singular
              inverse-square potentials, J. Funct. Anal., 254, (2008) (7),
              1864-1902.
\bibitem{judith} J. Vancostenoble and E. Zuazua,
Hardy inequalities, observability, and control for the wave
              and {S}chr\"odinger equations
               with singular potentials, SIAM J. Math. Anal., 41, 2009,
    4, 1508-1532.
\bibitem{vazzua} J. L. V\'{a}zquez and E. Zuazua, The Hardy inequality and the
asymptotic behaviour of the heat equation with an inverse-square
potential, J. Funct. Anal. 173 (2000), no. 1, 103-153.
\bibitem{hiddenenergy} J. L. V\'{a}zquez and N. B. Zographopoulos,
  Functional aspects of the Hardy inequality.
Appearance of a hidden energy, http://arxiv.org/abs/1102.5661.
\end{thebibliography}
\end{document}